\documentclass[
fleqn,oneside]{elsart}

\usepackage{amsmath}
\usepackage{amsfonts}
\usepackage{amssymb}
\usepackage{float}
\usepackage{natbib}
\usepackage{color}
\usepackage[]{hyperref}
\usepackage{subeqnarray}

\usepackage{anysize}

\def\bE{{\mathbf E}}

\def\bT{{\mathbf T}}
\def\cD{{\mathcal D}} 

\def\cT{{\mathcal T}}

\newcommand{\us}[1]{\usebox{#1}}
\newsavebox{\limn}
\sbox{\limn}{\large $\lim \limits_{n\mapsto\infty}$}

\newcounter{wspolnynrm}[section]
\newcounter{wspolnynrml}[section]

\newtheorem{Theorem}[wspolnynrm]{\textbf{Theorem}}
\newtheorem{Proposition}[wspolnynrm]{\textbf{Proposition}}
\newtheorem{Lemma}[wspolnynrml]{\textbf{Lemma}}
\newtheorem{Corollary}[wspolnynrml]{\textbf{Corollary}}
\newtheorem{remark}{\textbf{Remark}}
\journal{Stochastics}
\company{def}

\def\rightbox{\protect\vspace*{-2ex}

\begin{flushright}\(\blacksquare\)\end{flushright}}

\newenvironment{Proof}{{\bf Proof.}\hspace{1mm}}{\rightbox}

\begin{document}

\begin{frontmatter}
\title{The duration problem with multiple exchanges\thanksref{label1}}
\thanks[label1]{april 25, 2005}


\author{Charles E.M. Pearce\thanksref{label3}}
\address[label3]{The University of Adelaide, School of Mathematical Sciences,
 Adelaide, Australia SA 5005}
\ead{charles.pearce@adelaide.edu.au}
\ead[url]{http://www.maths.adelaide.edu.au/applied/staff/cpearce.html}

\author{Krzysztof Szajowski\corauthref{cor1}\thanksref{label2}}
\ead{Krzysztof.Szajowski@pwr.wroc.pl}
\ead[url]{http://neyman.im.pwr.wroc.pl/\~{}szajow}
\address[label2]{Wroc\l{}aw University of Tech., Institute of
Mathematics and Computer Science, Wybrze\.{z}e Wyspia\'{n}skiego 27, 50-370 Wroc\l{}aw, Poland}
\corauth[cor1]{Corresponding author}
\author{Mitsushi Tamaki\thanksref{label4}}
\ead{tamaki@vega.aichi-u.ac.jp}
\ead[url]{http://leo.aichi-u.ac.jp/\~{}~tamaki}
\address[label4]{Aichi Univ., Nagoya Campus: 370 Kurozasa, Miyoshi, Nishikamo, Aichi 470-02, Japan}

\date{ 21st of September 2006 }

\maketitle

\begin{abstract}
We treat a version of the multiple-choice secretary problem called the
multiple-choice duration problem, in which the objective is to maximize the time
of possession of relatively best objects. It is shown that, for the $m$--choice
duration problem, there exists a sequence $(s_1,s_2,\ldots,s_m)$ of critical
numbers such that, whenever there remain $k$ choices yet to be made, then the
optimal strategy immediately selects a relatively best object if it appears at
or after time $s_k$ ($1\leq k\leq m$). We also exhibit an equivalence between
the duration problem and the classical best-choice secretary problem. A simple
recursive formula is given for calculating the critical numbers when the number
of objects tends to infinity. Extensions are made to models involving an
acquisition or replacement cost.
\end{abstract}

\end{frontmatter}

\vskip-3cm
\renewcommand{\baselinestretch}{1.2}
\baselineskip 1.35pc

\section{Introduction and summary}
\cite{ferhartam92:own} were the first to consider a sequential and selection problem referred to 
as the duration problem, a variation of the classical secretary problem as treated by 
~\cite{gilmos66} and others (see ~\cite{fer89:who} and ~\cite{sam91:secretary} for a history and 
review of the secretary problem). The basic form of the duration problem in the no--information 
setting can be described as follows: a set of $n$ rankable objects appears one at a time in random 
order with all $n!$ permutations equally likely. As each objects appears, we decide either to 
select or reject it based on the relative ranks of the objects. The payoff is the length of time 
we are in possession of a relatively best object that have appeared to date. Thus we will select 
only a relatively best object, receiving unit payoff as we do so and an additional unit for each 
new observation, as long as the selected object remains relatively best.

Though ~\cite{ferhartam92:own} considered the various models extensively, they
confined themselves to the study of the one-choice problem. We consider here, as
a natural generalization of their study, the no-information multiple-choice
duration problem and its modifications. Preliminary results are included to
~\cite{tampeasza98:duration}. The multiple--choice duration problem is
reformulated as the multiple optimal stopping problem, which  has been treated
by many authors. The double--stopping problem was posed by ~\cite{hag67} and for
discrete--time Markov processes has been considered by ~\cite{eid79},
~\cite{nik79,nik98:multstopping} and ~\cite{sta85}.

For the $m$--choice duration problem, we choose at most $m$ objects sequentially, and receive
unit payoff at each time point as long as the last chosen object remains a {\it candidate},
that is, a relatively best object. Only candidates can be chosen, the objective being to maximize
the expected total payoff. More formally, this problem can be described as follows: let $\cT_n(i)$ 
denote the arrival time of the first candidate after time $i$ if there is one, and $n+1$ if there 
is none, so $\cT_n(i):\Omega\rightarrow \{i+1,\ldots,n+1\}$. Then
\begin{equation*}
\cD_n(i)\equiv \cT_n(i)-i
\end{equation*}
is the duration of the candidate selected at time $i$ and the objective is to find a stopping 
vector $(\tau_1^*,\tau_2^*,\ldots,\tau_m^*)$ such that
$$
 \bE\left[ \sum\limits_{i=1}^m \cD_n(\tau_i^*)\right]
 =\sup\limits_{(\tau_1^*,\tau_2^*,\ldots,\tau_m^*)\in C_m}
  \bE\left[ \sum\limits_{i=1}^m \cD_n(\tau_i)\right].
$$
Here $\tau_i$ ($1\leq i\leq m$) denotes the stopping time related to the $i$-th
choice and $C_m$ is the set of all possible vectors $(\tau_1,\tau_2,\ldots,\tau_m)$.

A generalization of this problem is considered in Section \ref{sec2BT}, in which we allow the 
number $M$ of objects presented to be a random variable. We show for the $m$--choice duration 
problem that, subject to a condition on the distribution of $M$,  there exists a nonincreasing 
sequence $(s_1,s_2,\ldots,s_m)$ of critical positive integers such that, whenever there remain $k$ 
choices to be made, the optimal strategy immediately selects a candidate if it appears at or after 
time $s_k$ ($1\leq k\leq m$). That is, the optimal strategy 
has a threshold form.
In Section \ref{sec2BT} we show that this condition is satisfied in the two particular cases when 
$M\sim M_s(n)$ (the distribution with $P\{ M=n\} =1$ \emph{i.e.} the degenerate distribution) and 
$M\sim M_u(n)$ (the distribution with $P\{ M=i\} =1/n\quad (i=1,\ldots ,n $ \emph{i.e.} the discrete uniform distribution).

In Section \ref{sec2BTA} we investigate the asymptotics for $n\rightarrow\infty$ in the case 
when $M$ has the degenrate distribution $P\{ M=n\} =1$. The ratio $s_k/n$ converges to a limit 
$s_k^*\in(0, 1)$. A recursive formula for calculating $s_k^*$ in terms of $s_1^*,s_2^*,\ldots, 
s_{k-1}^*$ is given by 
\begin{equation}
 s_k^*=\exp \left[-\left\{1+\sqrt{1-2\sum\limits_{i=1}^{k-1}
      \frac{[(k-i+2)+(k-i+1)\ln s_i^*]}{(k-i+2)!}(\ln s_i^*)^{k-i+1}}\right\}\right].
\label{row:1}
\end{equation}
Throughout the paper the empty sum is taken as zero, so the above formula is valid for all $k\geq 
1$. We show also that, as $n\rightarrow\infty$, the expected proportional payoff, that is, the 
expected maximum payoff per unit time, is given by  -$\sum\limits_{k=1}^m s_k^*\ln s_k^*$.

The classical best-choice secretary problem (\emph{BCSP}) is concerned with maximizing the 
probability of choosing the best object. ~\cite{sam91:secretary} and ~\cite{ferhartam92:own}
pointed out that, for the one-choice problem, the duration problem with a known number $n$
of objects is equivalent to the \emph{BCSP} with an unknown number of objects having a uniform 
distribution on $\{1,2,\ldots,n\}$. This was first studied by ~\cite{preson72:eng}. See also  
~\cite{pet83} and ~\cite{leh93:unknown} for the problem with an unknown number of objects,  
~\cite{gilmos66}, ~\cite{sak78d} and ~\cite{pre89:choice} for the multiple--choice problem and 
~\cite{tam79:ola} for the formulation of the multiple-choice problem with an unknown number of 
objects and solution of the two-choice problem having a uniform prior of the actual number of 
objects.

We show in Section \ref{sec2BT} that this equivalence still holds for the multiple-choice problem. 
Recently ~\cite{gne05:objectives} established an equivalence between the various best-choice 
problems and the related duration problems in a greater generality (see also 
~\cite{gne04:planar}). ~\cite{ferhartam92:own} considered another type of problem, the
full-information duration problem, where the observations are the actual values of the objects
assumed to be independent and identically distributed (iid) from a known distribution and hence 
decisions are based on the actual values of the objects. They showed that the above equivalence 
between the best-choice problem and the duration problem holds also for the full-information 
problem. ~\cite{por87:full,por02:similar} consider the full-information best-choice problem with a 
random number of objects and ~\cite{maztam06:dur} and ~\cite{sam04:why} the limiting maximum 
proportional payoff for the full-information one-choice duration problem.

In Sections~\ref{sec3} and \ref{sec4}, the multiple-choice duration problem with $M\sim M_s(n)$ is 
generalized by introducing costs. In Section~\ref{sec3} a constant acquisition cost is incurred 
each time an object is chosen, while in Section~\ref{sec4} a constant replacement cost is incurred 
with the selection of any candidate other than the first. The objective in Sections~\ref{sec3} and 
\ref{sec4} is to maximize the expected net payoff. It can be shown that, under an appropriate cost 
condition, the optimal strategies have similar structure to that for the problem involving no cost.
In Sections \ref{sec3asym} and \ref{sec4asym} we investigate the respective associated asymptotics.

The multiple-choice duration problem with replacement and acquisition costs may be considered as a 
marriage and divorce problem, interpreting the replacement cost as alimony. Recently,  
\cite{searap97:exper,searap00:unknown} investigated the behaviour of the decision makers under 
circumstances similar to those of the best choice problem model. The discussion of the problem are 
the subject of papers by \cite{bea06:cardinal} and \cite{sza06:risk}. It seems important to 
construct a model with physical parameter to fit it to the empirical data. The considerations of sections \ref{sec3} and \ref{sec4} suggest one way forward.

\vskip-3cm
\section{\label{sec2}Multiple exchange and hold of the relatively best item}
In this section we address the optimal choice problem with an unknown (bounded)
number of objects. At most $n$ objects appear in turn before us. We have only
an {\it a priori} distribution $p_i=P\{ M=i\}$ on the actual number of objects,
where $\sum _{i=1}^np_i=1$. Without loss of generality, we may assume that
$p_n>0$, so $\pi _i>0$ for $1\leq i\leq n$. We set $\sum _{j=i}^np_i=\pi _i$. We
are allowed to make at most $m$ choices and wish to maximize the expected
duration of holding a relatively best object.

We  assume that all that can be observed are the relative ranks of the objects as they are 
presented. Thus if $X_i$ denotes the relative rank of the $i$-th object amongst those observed so 
far (a candidate if $X_i=1$), the sequentially-observed random variables are $X_1,X_2,\ldots,X_n$. 
It is well known that under the assumption that the objects are in random order with all $n!$ 
permutations equally likely, we have
\begin{description}
        \item[(a)] the $X_i$ are independent random variables and
        \item[(b)] P\{$X_i=j$\}=$1/i$, for $1\leq j\leq i$, $1\leq i\leq n$.
\end{description}

We formulate the $m$-choice duration problem as a Markovian decision-process
model. First we condition on $M=\ell \geq i$. Since decisions about selection or
rejection occur only when a candidate appears, we describe the state of the
decision process as $(i,\,k)$, $1\leq i\leq \ell$, $1\leq k\leq m$, if the
$i$-th object is a candidate and there remain $k$ more choices to be made. For
the above process to be a Markov chain, we must further introduce an additional
absorbing state $(\ell + 1, k,e)$ for the situation where the last object is
presented at time $\ell$ and is not a candidate, with $k$ choices left
($1\leq k\leq m$). When it leaves $(i,k)$, the process moves to a state $(j,k-1)$ if the $i$-th 
object is selected. Otherwise it moves to a state $(j,k)$ or $(\ell +1,k,e)$.
By (a) and (b), the distribution of $j$ is given by
\begin{equation}
 p(i,\,j)=\begin{cases} \frac{i}{j(j-1)}, &j=i+1,\ldots,\ell \\
                \frac{i}{\ell},&j=\ell +1
          \end{cases}  .
\label{row:2}
\end{equation}
We now remove the conditioning on $M$. The probability of a next candidate
appearing at time $j$ ($i<j\leq n$) is given by
\begin{equation}\label{condPROB}
p(i,j)=\frac{i}{j(j-1)}\cdot \frac{\pi _j}{\pi _i}.
\end{equation}
The corresponding probability of transition from $(i,k)$ to a state $(j,k)$,
$(j,k-1)$ or $(j,k,e)$ is
\begin{equation}\label{tranPROB}
q(i,j)=\frac{i}{j(j-1)}\cdot \frac{\pi _j}{\pi _i}+\frac{i}{j-1}\cdot
\frac{p_{j-1}}{\pi _i}.
\end{equation}
In accordance with our convention about empty sums, we shall interpret
$\pi_{n+1}$ as zero. 

Let $f:\{1,2,\ldots,n+1\}\rightarrow \Re$ be the payoff function. Define
\begin{equation}\label{opT}
\bT f(i)=\sum\limits_{j=i+1}^{n+1}f(j)p(i,\,j)
\end{equation}
the expectation operator with respect of the probability distribution (\ref{condPROB})
and 
\begin{equation}\label{opTE}
\bT^{e} f(i)=\sum\limits_{j=i+1}^{n+1}f(j)q(i,\,j)
\end{equation}
the expectation operator with respect of the probability distribution (\ref{tranPROB}).

If the $i$-th object observed is a candidate, the period for which it remains a candidate has mean
\begin{equation}
U_i^{(1)} =\bT^{e}\cD_n(i)= \sum _{j=i+1}^{n+1}(j-i)q(i,j)=\frac{i}{\pi _i}\sum _{j=i}^n
\frac{\pi _j}{j}=\frac{i}{\pi _i}(L_n-L_{i-1}),
\label{dtilde}
\end{equation}
where
\begin{equation}\label{Lsum}
L_j=\sum _{\ell =1}^j\frac{\pi _{\ell}}{\ell}\ .
\end{equation}
When $\pi_j=1$, for $j=1,2,\ldots,n$ then the sequence $L_j$ given in (\ref{Lsum}) is denoted by 
$H_j$ (see section \ref{sec5} for further details and section \ref{sec2BTA} for properties given 
by formulae (\ref{eq:G1H})). The contribution to the expected occupancy time from a further 
candidate, if any, is
\begin{equation}
\widetilde{V}_i^{(1)}=\bT[U_{\cT_n(i)}^{(1)}]=\sum _{j=i+1}^np(i,j)U_j^{(1)}
=\frac{i}{\pi _i}\sum _{j=i+1}^n\frac{1}{j-1}\left( L_n-L_{j-1}\right) .
\label{eq:A16}
\end{equation}

When $n\leq m$, the optimal strategy is easily seen to be to select the
candidates successively as they appear. Thus we assume $n>m$. Before proceeding
to investigate the optimal strategy, we introduce some notation. Suppose we
start in state $(i,k)$. We denote by $U_i^{(k)}$, $V_i^{(k)}$ the expected total
possession time when we select or reject respectively the $i$-th object and then
proceed in an optimal manner. We also denote by $W_i^{(k)}$ the expected total
possession time under an optimal strategy starting from state $(i,k)$ ($1\leq
i\leq n$, $1\leq k\leq m$). The Bellman principle of optimality yields
for $1\leq i\leq n$ and $1\leq k\leq m$ that
\begin{eqnarray}
W_i^{(k)}&=&\max \{ U_i^{(k)},V_i^{(k)}\} ,
\label{row:3} \\
U_i^{(k)}&=&U_i^{(1)}+\sum _{j=i+1}^np(i,j) W_j^{(k-1)},
\label{row:4} \\
V_i^{(k)}&=&\sum _{j=i+1}^np(i,j)W_j^{(k)}.
\label{row:5}
\end{eqnarray}
Equations (\ref{row:3})--(\ref{row:5}) together with the boundary conditions
$W_i^{(0)}=0$ for $1\leq i\leq n+1$ and $W_{n+1}^{(k)}=0$ for $1\leq k\leq m$
can be solved recursively to yield the optimal strategy and the optimal value
$W_1^{(m)}$. When $n\leq m$ selecting the candidates successively as they
appear is the optimal strategy and achieves total possession time $n$ for
relatively best objects. Thus we may assume $n>m$ without loss of generality.

\subsection{\label{sec2tex}The main auxiliary theorem}
In the sequel, in the construction of optimal solutions, the properties of the sequence of differences between some payoffs and expected payoffs will be analysed. Important properties of such sequences are gathered in the following theorem. 

\begin{Theorem}\label{thm11}  Suppose that $s_1$ and $N$ are integers with $1\leq s_1\leq N$. Suppose further that $(G_i^{(1)})_{i=1}^{N}$ is such that
\begin{description}
\item[$(C1)$] if $G_i^{(1)} \geq 0$ with $1\leq i<N$, then $G_{i+1}^{(1)}\geq 0$;
\item[$(C2)$] if $1\leq i<s_1$, then $G_{i+1}^{(1)}\geq G_i^{(1)}$;
\item[$(C3)$] $s_1$ is determined by
$$  s_1=\min\{i:G_i^{(1)}\geq 0\}.$$
\end{description}

Then if $1\leq m\leq N$ there exists a sequence $(s_1,s_2,\ldots,s_m)$ of positive integers
such that if, for $1\leq i\leq N$, we define
 \begin{equation}
G_i^{(k)}=G_i^{(1)}+ \sum\limits_{j=\max(i+1, s_{k-1})}^{N}
\frac{1}{j-1}G_j^{(k-1)} \hspace{1cm} (k \geq 2)
 \label{row:7}
 \end{equation}
recursively for $1<k\leq m$, then we have the properties:
\begin{description}
\item[$(P1)_k$] if $G_i^{(k)} \geq 0$ with $1\leq i<N$, then $G_{i+1}^{(k)}\geq 0$;
\item[$(P2)_k$] if $1\leq k<m$, then $G_i^{(k+1)} \geq G_i^{(k)}$;
\item[$(P3)_k$] if $1\leq i<s_k$, then $G_{i+1}^{(k)}\geq G_i^{(k)}$;
\item[$(P4)_k$] $s_k$ is determined by
 \begin{equation}
  s_k=\min\{i:G_i^{(k)}\geq 0\};
  \label{row:6}
 \end{equation}
\item[$(P5)_k$] $(s_1,\ldots ,s_{k+1})$ is nonincreasing for $1\leq k<m$.
\end{description}
\end{Theorem}

\begin{Proof}
We shall employ induction on $k$. From $(C1)$--$(C3)$ we have
$(P1)_1$, $(P3)_1$ and $(P4)_1$
and that $G_i^{(2)}$ is well-defined for $1\leq i\leq N$, with
\begin{equation}
G_i^{(2)}=G_i^{(1)}+ \sum\limits_{j=\max(i+1, s_1)}^{N}
\frac{1}{j-1}G_j^{(1)}.
\label{eq:**}
\end{equation}
By $(P1)_1$ and $(P4)_1$ the summand in
(\ref{eq:**}) is nonnegative so that (\ref{eq:**}) yields
$(P2)_1$. From this and $(P4)_1$ we deduce $(P5)_1$. Thus we have a basis for
the induction.

 For the inductive step, suppose that for some $k$ with $1\leq k<m$
$(P1)_{\ell}$--$(P5)_{\ell}$ hold for $1\leq \ell \leq k$.
By definition we have for $1\leq i\leq N$ that
$$
G_i^{(k+2)}-G_i^{(k+1)}=\sum\limits_{j=\max(i+1, s_{k+1})}^{N}
\frac{1}{j-1}G_j^{(k+1)}-\sum\limits_{j=\max(i+1, s_k)}^{N}
\frac{1}{j-1}G_j^{(k)}                         .
$$
 From $(P1)_{k}$ and $(P4)_{k}$ the second summand is nonnegative. Since
$s_k\geq s_{k+1}$, we thus have
$$
G_i^{(k+2)}-G_i^{(k+1)}\geq \sum\limits_{j=\max(i+1, s_{k+1})}^{N}
\frac{1}{j-1}\left[ G_j^{(k+1)}- G_j^{(k)}\right] .
$$
Hence $(P2)_{k}$ implies that
$$G_i^{(k+2)}-G_i^{(k+1)}\geq 0$$
and we have $(P2)_{k+1}$. Thus $s_{k+2}$ is well-defined by (\ref{row:6})
and $(P4)_{k+1}$ applies. Also $s_{k+2}\leq s_{k+1}$, so $(P5)_{k+1}$ holds.

Since $s_{k}\leq s_1$, $(P3)_1$ implies $G_{i+1}^{(1)}\geq G_i^{(1)}$ for
$1\leq i<s_{k}-1$. For any such values of $i$,
$$G_{i+1}^{(k+1)}=G_{i+1}^{(1)}+ \sum\limits_{j=s_{k}}^{N}
\frac{1}{j-1}G_j^{(k)},$$
$$G_{i}^{(k+1)}=G_{i}^{(1)}+ \sum\limits_{j=s_{k}}^{N}
\frac{1}{j-1}G_j^{(k)}.$$
Thus
\begin{equation}
G_{i+1}^{(k+1)}\geq G_{i}^{(k+1)}\quad \mbox{ for $1\leq i<s_{k}-1$}.
\label{eq:AA}
\end{equation}
 Also by
(\ref{row:6}) we have for $i=s_{k+1}-1$ that $G_i^{(k+1)}<0$ and
$G_{i+1}^{(k+1)}\geq 0$. Hence $(P3)_{k+1}$ applies.

This leaves $(P1)_{k+1}$.
Since $G_{s_{k+1}}^{(k+1)}\geq 0$, we have by (\ref{eq:AA}) that
$G_{i}^{(k+1)}\geq 0$
for $s_{k+1}\leq i< s_k$. Further by $(P2)_{k}$ we have
$G_{i}^{(k+1)}\geq 0$ for $i\geq s_k$. Hence $G_{i}^{(k+1)}\geq 0$
 for $i\geq s_{k+1}$ and $G_i^{(k+1)}<0$ for $i<s_{k+1}$, giving $(P1)_{k+1}$.
This completes the inductive step.
\end{Proof}

\subsection{\label{sec2basic}The basic theorem}
We complement relations (\ref{row:3})--(\ref{row:5}) with $\widetilde{V}_i^{(k)}$, the expected 
total possession time if a candidate at time $i$ is rejected and the next candidate (if any) 
accepted, with optimal choices following such acceptance.

\begin{Theorem}\label{thm21} For $1\leq i\leq n$, let $G_i^{(1)}$ be given by (\ref{eq:UVG}) 
below with $k=1$. Suppose that for $N=n$ and some integer $s_1$, conditions $(C1)$--$(C3)$
of Theorem \ref{thm11} are satisfied. Then for the $m$-choice duration problem with $k$
($1\leq k\leq m$) choices still to be made, the optimal strategy selects
immediately the first candidate, if any,  to appear at or after time $s_k$ of Theorem \ref{thm11}.
\end{Theorem}

\begin{Proof} Since conditions $(C1)$--$(C3)$
of Theorem \ref{thm11} are satisfied, the conclusions of Theorem \ref{thm11} hold.
We proceed inductively, establishing the following:
\begin{description}
\item[$(Q1)_k$] the optimal strategy when there are $k$ choices still to be
made to to select the first candidate to appear at or after time $s_k$;
\item[$(Q2)_k$] for $1\leq i\leq n$
\begin{equation}
\frac{\pi _i}{i}\left[ U_i^{(k)}-\widetilde{V}_i^{(k)}\right] =G_i^{(k)}.
\label{eq:UVG}
\end{equation}
\end{description}

 For a basis, consider the one--choice duration problem. By definition
$(Q2)_1$ is given. We should select a
candidate observed at time $i$ in preference to the next candidate (if any) if
$U_i^{(1)}\geq \widetilde{V}_i^{(1)}$,
that is, if
\begin{equation}
G_i^{(1)}\geq 0.
\label{eq:UVG1}
\end{equation}
By $(P4)_1$ of Theorem \ref{thm11}, this condition cannot be satisfied if $i<s_k$.
Also by $(P1)_1$ and $(P4)_1$, the first candidate at or after time $s_1$
satisfies (\ref{eq:UVG1}) and choice of this candidate is strictly preferable
to choice of any candidate subsequent to the second candidate after $s_k$.
Thus $(Q1)_1$ holds and $k=1$ provides a basis for induction.

 For the inductive step, suppose $(Q1)_{\ell}$ and $(Q2)_{\ell}$ to be true
for $\ell =1,\ldots ,k-1$ for some $k$ with $2\leq k\leq m$. Then
\begin{equation}
  W_j^{(k-1)} =
    \begin{cases}
      V_j^{(k-1)}, & j <s_{k-1}\\
      U_j^{(k-1)}, & j \geq s_{k-1}
    \end{cases}
\label{eq:WUV}
\end{equation}
and
\begin{equation}
  V_j^{(k-1)} = \widetilde{V}_j^{(k-1)},\ \ j \geq s_{k-1} - 1.
\label{eq:VV*}
\end{equation}
By the Bellman principle of optimality we have
\begin{equation}
\widetilde{V}_i^{(k)}=\sum _{j=i+1}^np(i,j)\left[ U_j^{(1)}+V_j^{(k-1)}
\right]
=\widetilde{V}_i^{(1)}+\sum _{j=i+1}^np(i,j)V_j^{(k-1)}
\label{eq:VVV}
\end{equation}
by (\ref{eq:A16}).

By (\ref{eq:WUV}), subtraction of (\ref{eq:VVV}) from (\ref{row:4}) yields
\begin{eqnarray*}
U_i^{(k)}-\widetilde{V}_i^{(k)}&=&U_i^{(1)}-\widetilde{V}_i^{(1)}
+\sum _{j=\max (i+1,s_{k-1}) }^np(i,j)\left[ U_j^{(k-1)}-V_j
^{(k-1)}\right] \\
&=&U_i^{(1)}-\widetilde{V}_i^{(1)}
+\sum _{j=\max (i+1,s_{k-1}) }^np(i,j)\left[ U_j^{(k-1)}-\widetilde{V}_j
^{(k-1)}\right] \\
\end{eqnarray*}
(by (\ref{eq:VV*})). Thus the inductive assumption provides
$$\frac{\pi _i}{i}\left[ U_i^{(k)}-\widetilde{V}_i^{(k)}\right] =
G_i^{(1)}+\sum _{j=\max ( i+1,s_{k-1}) }^n\frac{1}{j-1}G_j^{(k-1)}.$$
The recursive definition of the functions $G_i^{(k)}$ leads to
$$\frac{\pi _i}{i}\left[ U_i^{(k)}-\widetilde{V}_i^{(k)}\right] = G_i^{(k)}$$
and we have established $(Q2)_{k}$. The argument leading from $(Q2)_{k}$
to $(Q1)_{k}$ follows that leading from $(Q2)_{1}$
to $(Q1)_{1}$ and the inductive step is
complete.
\end{Proof}

\subsection{\label{sec2BT}Applications of the basic theorem}
We can apply Theorem \ref{thm21} whenever we can verify conditions $(C1)$--$(C3)$.
To this end, we note that by (\ref{dtilde}) and (\ref{eq:A16})
\begin{equation}
G_i^{(1)}=L_n-L_{i-1}-\sum _{j=i+1}^n\frac{1}{j-1}(L_n-L_{j-1})
\label{eq:G1}
\end{equation}
for $1\leq i\leq n$, whence we derive that
\begin{equation}
G_{i+1}^{(1)}\gtreqless G_{i}^{(1)}\quad \mbox{ according as }\quad
L_n \gtreqless L_i+\pi _i\ =:\phi _i
\label{eq:gel}
\end{equation}
for $1\leq i<n$.

\begin{Proposition}\label{prop41}
For $1\leq i\leq n$ put $\psi _i=(i+1)\pi _i$. Then a sufficient condition
for $(C1)$--$(C3)$ to hold is that there should exist an integer $i_0$ with
$1\leq i_0\leq n$ such that
\begin{equation}
\psi _1\leq \psi _2 \leq \ldots \leq \psi _{i_0}\quad \mbox{ and }\quad
\psi _{i_0}\geq \psi _{i_0+1}\geq \ldots \geq \psi _n.
\label{eq:psi}
\end{equation}
\end{Proposition}
\begin{Proof}
Since $L_n<\phi _n$ and $G_n^{(1)}=\pi _n/n>0$, we have readily from (\ref{eq:gel})
that if $1\leq i_0\leq n$ with
$$\phi _1\leq \phi _2 \leq \ldots \leq \phi _{i_0}\quad \mbox{ and }\quad
\phi _{i_0}\geq \phi _{i_0+1}\geq \ldots \geq \phi _n,$$
then $G_i^{(1)}$ satisfies $(C1)$--$(C3)$.  Also
$\phi _{i+1} \gtreqless \phi _i\quad \mbox{ according as }\quad
\psi _{i+1} \gtreqless \psi _i$.
The stated result follows.
\end{Proof}

As corollaries we consider the two special choices $M\sim M_s(n)$ (when $M=n$
with probability one) and $M\sim M_u(n)$ (when $M$ is uniformly distributed
on $\{ 1,2,\ldots ,n\} $).

\begin{Corollary}\label{corol11}
For $1\leq m<n$, the optimal strategy in the $m$--choice
duration problem with $M\sim M_s(n)$ is given by the conclusion of Theorem \ref{thm21}.
\end{Corollary}

\begin{Proof} Here $\pi _i=1$ for $1\leq i\leq n$, so
(\ref{eq:psi}) holds with $i_0=n$. The stated result follows from Theorem~\ref{thm21}
and Proposition~\ref{prop41}.
\end{Proof}

\begin{Corollary}\label{corol12}
For $1\leq m<n$ and $M\sim M_u(n)$, the optimal strategy in the $m$--choice
duration problem is given by the conclusion of Theorem \ref{thm21}.
\end{Corollary}

\begin{Proof} Here $\pi _i=(n-i+1)/n$ for $1\leq i\leq n$, so
$\psi _i=(i+1)(n-i+1)/n$ and
$$
\psi _{i+1}\gtreqless \psi _{i}\quad \mbox{ according as }\quad n \gtreqless 2i+1.
$$
Hence (\ref{eq:psi}) holds with $i_0=\max \{ 1,\lfloor (n-1)/2\rfloor \}$ and
 the stated result follows from Theorem \ref{thm21} and Proposition \ref{prop41}.
\end{Proof}

 For $m\geq 1$, the expected payoff for the $m$--choice problem with
$M\sim M_s(n)$ is
$q_m\equiv W_1^{(m)}$. We have
\begin{equation}
q_m=\widetilde{V}_{s_m-1}^{(m)}
=\frac{s_m-1}{n}\sum _{j=s_m}^n\frac{1}{j-1}\sum_{t=j}^n\frac{1}{t}+
\sum_{j=s_m}^n\frac{s_m-1}{j(j-1)})V_j^{(m-1)}.
\label{row:14}
\end{equation}
The $V_i^{(m-1)}$ ($m\geq 2$) may be calculated recursively from
\begin{equation}
V_i^{(m-1)}=\begin{cases}
             q_{m-1}, & i<s_{m-1}-1\\
             \widetilde{V}_{i}^{(1)}+\sum _{j=i+1}^np(i,j)V_j^{(m-2)}, &
\mbox{ $ i\geq s_{m-1}-1$}
             \end{cases}
\label{row:15}
\end{equation}
with the interpretation that $V_i^{(0)}\equiv 0$.

In the three subsequent sections on asymptotics it is convenient to scale mean durations by 
dividing by $n$ so as to work in terms of the average possession times per unit time. We shall set 
$\overline{U}_i^{(k)}=U_i^{(k)}/n$, etc. This leaves optimal strategies unaffected. As a
prelude to this, we consider one further application.

We compare two differently  formulated multiple optimal stopping problems.
The first is the $m$-choice problem of Corollary \ref{corol11} and the second the $m$-choice
 best-secretary problem with an unknown number of objects having distribution $M_u(n)$. We show 
that these problems have the same solution in the sense that the optimal strategies and expected 
payoffs are the same. In the latter problem we win if the last chosen object is
best overall. The objective is to maximize the winning probability.

\cite{tam79:ola} formulated this problem as a Markovian decision process model
and solved explicitly the two-choice problem with a uniform prior on $M$. We
describe the state of the process be described as $(i,k)$  ($1 \leq i \leq n$,
$1 \leq k \leq m$) if the $i$-th object is a candidate and there remain $k$
choices to be made. We denote by $u_i^{(k)}(v_i^{(k)})$ the winning probability
when we select (reject) the $i$-th object and then  continue optimally from
state $(i,k)$. If we let
\begin{equation}\label{row:31}
  w_i^{(k)} = \max\left\{u_i^{(k)},v_i^{(k)}\right\},\quad \mbox{$1 \leq i$} \leq n
\end{equation}
the principle of optimality yields
\begin{eqnarray}\label{row:32}
  u_i^{(k)} & = & \sum\limits_{j=i}^n \frac{i}{j}\frac{p_j}{\pi_i} +
                  \sum\limits_{j=i+1}^n \frac{i}{j(j-1)} \frac{\pi_j}{\pi_i}
                  w_j^{(k-1)},\\
                \label{row:33}
  v_i^{(k)} & = & \sum\limits_{j=i+1}^n \frac{i}{j(j-1)}\frac{\pi_j}{\pi_i} w_j^{(k)}.
\end{eqnarray}
When $M\sim M_u(n)$, we have $\pi _i=(n-i+1)/n$ for $i=1,\ldots ,n$.
If we set
\begin{equation}
\overline{U}_i^{(k)} = \frac{n-i+1}{n}  u_i^{(k)},\quad \overline{V}_i^{(k)} =
\frac{n-i+1}{n} v_i^{(k)}\quad \mbox{ and }\quad \overline{W}_i^{(k)}
=\frac{n-i+1}{n} w_i^{(k)},
\label{row:ch}
\end{equation}
(\ref{row:31})--(\ref{row:33}) are transformed respectively into
\begin{eqnarray*}
\overline{W}_i^{(k)}&=&\max \{ \overline{U}_i^{(k)},\overline{V}_i^{(k)}\} ,\\
\overline{U}_i^{(k)}&=&\frac{i}{n}\sum _{j=i}^n\frac{1}{j}+\sum
_{j=i+1}^n\frac{i}{j(j-1)} \overline{W}_j^{(k-1)},\\
 \overline{V}_i^{(k)}&=&\sum _{j=i+1}^n\frac{i}{j(j-1)}\overline{W}_j^{(k)}   .
\end{eqnarray*}
These are (\ref{row:3})--(\ref{row:5}) for the scaled version of the process
of Corollary \ref{corol11}, with the correct normalized  value for $\overline{U}_i^{(1)}$.
 Because of the common multiplicative factors in (\ref{row:ch}), $w_i^{(k)}=
u_i^{(k)}$ if and only if $W_i^{(k)}= U_i^{(k)}$. Thus optimal choices are the
same in the two processes. Since $W_1^{(k)}=w_1^{(k)}$, the two also share
their optimal payoff value. Thus we have established the following result.

\begin{Theorem}\label{thm31} 
The optimal-choice strategy and expected payoff is the same for the $m$-choice
versions of
\begin{description}
\item[$(i)$] the best-choice secretary problem with unknown number of objects
distributed uniformly on $\{ 1,2,\ldots ,n\}$;
\item[$(ii)$] the duration problem with possession times of relatively best
objects scaled by division by $n$.
\end{description}
\end{Theorem}

\vskip-2.5cm
\subsection{\label{sec2BTA}Asymptotics for the basic problem with the degenerated distribution of objects}
It is of interest to investigate the asymptotic behaviour of $s_k/n$ ($1\leq k\leq m$) and $q_m/n$ 
as $n\to \infty$. To do this, we observe that the sums in the formula of Section 2 are Riemann 
sums. With $M\sim M_s(n)$, (\ref{eq:G1}) becomes
\begin{equation}
G_i^{(1)}=H_n-H_{i-1}-\sum _{j=i+1}^n\frac{1}{j-1}(H_n-H_{j-1})
\label{eq:G1H}
\end{equation}
for $1\leq i\leq n$, where $H_{\ell} =\sum_{i=\ell}^n1/{\ell}$ for $\ell \geq 1$ and $H_0=0$ (see also (\ref{Lsum})).

For $m=1$ with $i/n\rightarrow x$ as $n\rightarrow\infty$,
the Riemann sum given by $G_i^{(1)}$ converges to the integral
\begin{equation}
 G^{(1)}(x)=\int_x^1\frac{dy}{y}-\int_x^1\frac{dy}{y}\int_y^1\frac{dz}{z}=
  -\frac{(2+\ln x)\ln x}{2}.
\label{row:16}
\end{equation}
 From (\ref{row:6}),
\begin{equation}
 s_1^*=\us{\limn}s_1/n=e^{-2}
\label{eq:s1}
\end{equation}
 is obtained
as the unique root $x\in(0,\,1)$ of $G^{(1)}(x)=0$.

More generally (\ref{row:7}) leads to functions $G^{(k)}(x)$ ($0<x<1$)
defined recursively by
\begin{equation}
 G^{(k)}(x)=G^{(1)}(x)+\int_{\max(x,\, s_{k-1}^*)}^1 \frac{1}{y}G^{(k-1)}(y)dy,
 \hspace{1cm} k\geq 2
\label{row:18}
\end{equation}
with $G_i^{(k)}$ a Riemann
approximation to $G^{(k)}(x)$ if $i/n\rightarrow x$ as
$n\rightarrow\infty$\,.

Correspondingly  $s_k^*:=\us{\limn}s_k/n$ exists and
may be obtained for $k\geq 2$ as the unique root
$x\in(0,\,s_{k-1}^*)$ of
\begin{equation}
 G^{(k)}(x)=0.
\label{row:17}
\end{equation}

 From (\ref{row:18}) and (\ref{row:17}), $s_k^*$ is a root of
$$
 G^{(1)}(x)=-\int_{s_{k-1}^*}^1 \frac{1}{y}G^{(k-1)}(y)dy,
$$
or equivalently, from (\ref{row:16}),
\begin{equation}
 s_k^*=\exp\left\{-\left(1+\sqrt{1+2\int_{s_{k-1}^*}^1\frac{G^{(k-1)}(y)}{y}dy}
 \right)\right\}.
\label{row:19}
\end{equation}
To derive the tractable form (\ref{row:1}), we need
some lemmata.

\begin{Lemma}\label{lem21}
 For $k$ a positive integer, define
\begin{eqnarray}
A_{k,\,i}&=&\int_{s_{k-i}^*}^1\frac{(\ln x)^i}{x}\, G^{(k-i)}(x)dx,\mbox{
$0\leq i\leq k-1$,}\nonumber \\
&& \label{row:20} \\
a_{k,\,i}&=&\int_{s_{k-i}^*}^1\frac{(\ln x)^i}{x}\,  G^{(1)}(x)dx, \mbox{
$0\leq i\leq k-1$.} \nonumber
\end{eqnarray}
Then $A_{k,\,i}$ satisfies the recursion
\begin{equation}\label{row:21}
A_{k,\,i}=a_{k,\,i}+\frac{1}{i+1}\left[ A_{k,\,i+1}-(\ln
s_{k-i}^*)^{i+1} A_{k-i-1,\,0}\right] ,
\end{equation}
with $A_{k,\,k}=0$ ($k\geq 0$).
\end{Lemma}

\begin{Proof}
From (\ref{row:18})
$$
G^{(k-i)}(x)=
   \begin{cases}
      G^{(1)}(x)+A_{k-i-1,\,0},&\text{ $x< s_{k-i-1}^*$} \\
G^{(1)}(x)+\int_x^1\frac{1}{y}G^{(k-i-1)}(y)dy,&\text{  $x\geq s_{k-i-1}^*$}.
   \end{cases}
$$
Since $s_{k-i}^* \leq s_{k-i-1}^* $, we have
\small
\begin{eqnarray*}
A_{k,\,i} &=&\int_{s_{k-i}^*}^{s_{k-i-1}^*}\frac{(\ln
x)^i}{x}[G^{(1)}(x)+A_{k-i-1,\,0}]dx \\
\nonumber&&\mbox{}+ \int_{s_{k-i}^*}^1\frac{(\ln
x)^i}{x}\left[G^{(1)}(x)+\int_x^1\frac{1}{y}\,G^{(k-i-1)}(y)dy\right]dx\\
\nonumber& =&a_{k,\,i}+ A_{k-i-1,\,0}\int_{s_{k-i}^*}^{s_{k-i-1}^*}\frac{(\ln
x)^i}{x}dx
+\int_{s_{k-i}^*}^1\left[\int_{s_{k-i}^*}^y\frac{(\ln
x)^i}{x}dx\right]\frac{1}{y}\,G^{(k-i-1)}(y)dy.
\end{eqnarray*}

\normalsize

The second and third terms of the last line are
respectively
$$
\frac{A_{k-i-1,\,0}}{i+1}\left[ (\ln s_{k-i-1}^*)^{i+1}-
(\ln s_{k-i}^*)^{i+1}\right]
$$
and
\begin{eqnarray*}
&&\frac{1}{i+1}\int_{s_{k-i-1}^*}^1\left[ (\ln y)^{i+1}-(\ln
  s_{k-i-1}^*)^{i+1}\right] \frac{1}{y}\,G^{(k-i-1)}(y)dy\\
&&\quad\quad =\frac{1}{i+1}\left[ A_{k,\,i+1}-A_{k-i-1,\,0}\,
 (\ln s_{k-i-1}^*)^{i+1}\right] ,
\end{eqnarray*}
whence the desired result.
\end{Proof}

 For simplicity, set $A_k:=A_{k,\,0}$. Repeated use of
(\ref{row:21}) gives the following recursion for $A_k$.

\begin{Lemma}\label{lem22}
 For $k\geq 1$, $A_k$ satisfies the recursion
$$
 A_k=\sum\limits_{i=1}^k\left[\frac{a_{k,\,k-1}}{(k-i)!}-
   \frac{(\ln s_i^*)^{k-i+1}}{(k-i+1)!}A_{i-1}\right].
$$
\end{Lemma}

 For $k\geq1$, define $N_k$ by
\begin{equation}
 N_k=-(1+\sqrt{1+2A_{k-1}}).
\label{row:23}
\end{equation}
Then from (\ref{row:19})
\begin{equation}
 s_k^*=\exp \ N_k
\label{row:24}
\end{equation}
and we have the following lemma.

\begin{Lemma}\label{lem23}
 For $k\geq 1$, $N_k$ satisfies the recursion
\begin{equation}
 N_k=-\left[1+\sqrt{1-2\sum\limits_{i=1}^{k-1}
  \frac{\{(k-i+2)+(k-i+1)N_i\}(N_i)^{k-i+1}}{(k-i+2)!}}\right].
\label{row:25}
\end{equation}
\end{Lemma}

\begin{Proof}
 Straightforward calculation from
(\ref{row:20}) yields
$$  a_{k-1,\,k-i-1} = \frac{N_i^{k-i-1}}{k-i+1} +
    \frac{N_i^{k-i+2}}{2(k-i+2)}\quad \mbox{ and }\quad
  A_{i-1} = N_i + N_i^2/2,
$$
so from
Lemma~\ref{lem22}
\begin{eqnarray*}
  A_{k-1} & = & \sum\limits_{i=1}^{k-1} \left[
\frac{a_{k-1,\,k-i-1}}{(k-i-1)!} - \frac{(N_i)^{k-i}}{(k-1)!} A_{i-1}\right]\\
& = & -\sum\limits_{i=1}^{k-1}\frac{\{(k-i+2) +
(k-i+1)N_i\}(N_i)^{k-i+1}}{(k-i+2)!} \ .
\end{eqnarray*}
Combining this with (\ref{row:23}) completes the proof.
\end{Proof}

By (\ref{row:24}), recursion (\ref{row:1}) is an immediate consequence
of (\ref{row:25}). From (\ref{row:1}) we
have successively
\begin{eqnarray*}
  s_1^* & = & \exp\{-2\} \approx 0.1353,\\
  s_2^* & = & \exp\left\{
    -\left(1+\sqrt{\frac{7}{3}}\right) \right\} \approx
    0.0799,\\
  s_3^* & = & \exp\left\{
    -\left(1+\frac{1}{3}\sqrt{15+14\sqrt{\frac{7}{3}}}\right)
    \right\} \approx 0.0493,\\
  s_4^* & = & \exp\left\{
    -\left(1+\sqrt{\frac{31}{45} +
    \frac{2}{81}\left(15+14\sqrt{\frac{7}{3}}\right)^{3/2}}
    \right)\right\} \approx 0.0311.
\end{eqnarray*}
See Table \ref{tab1} for $s_5^*$ and $s_{10}^*$ ($c=0$).

\par We have the following lemma concerning the expected payoff.

\begin{Lemma}\label{lem24}
Let $q^*_m = \lim\limits_{n\rightarrow \infty} q_m/n$
for $m \geq 1$. Then
\begin{equation}
  \label{row:28}
  q_m^* = - \sum\limits_{k=1}^m s_k^* \ln s_k^*.
\end{equation}
\end{Lemma}

\begin{Proof}
For $m=1$, we have from (\ref{row:14})
that
\[
q_1^* = s_1^* \int\limits_{s_1^*}^{1} \frac{dy}{y}\int\limits_{y}^1 \frac{dz}{z}
= \frac{s_1^*}{2}(\ln s_1^*)^2.
\]
By (\ref{eq:s1}), this may also be written
\begin{equation}
q_1^* = -s_1^*\ln s_1^*.
\label{eq:qs}
\end{equation}
For $m \geq 2$, we have from (\ref{row:14}) and (\ref{row:15}) that
\begin{eqnarray}
  q_m^* & = & s_m^* \int\limits_{s_m^*}^{1} \frac{dy}{y}
    \int\limits_{y}^1 \frac{dz}{z} + s_m^* \int\limits_{s_m^*}^1
    \frac{1}{y^2} V^{(m-1)}(y) dy \nonumber\\
\label{row:29}&&\\
& = & \frac{s_m^*}{2}(\ln s_m^*)^2 + s_m^* \int\limits_{s_m^*}^1 \frac{1}{y^2}
V^{(m-1)}(y) dy\nonumber,
\end{eqnarray}
if $V^{(m-1)}(x)$ ($0 < x < 1,\ m \geq 2$) are defined recursively by
\[
  V^{(m-1)}(x) = \begin{cases}
    q_{m-1}^*, &\mbox{ $0 < x < s_{m-1}^*$,}\\
x \int\limits_{x}^{1} \frac{dy}{y}\int\limits_{y}^1 \frac{dz}{z} + x
\int\limits_{x}^1
    \frac{1}{y^2} V^{(m-2)}(y) dy, &\mbox{ $s_{m-1}^* \leq x < 1$}
  \end{cases}
\]
starting with $V^{(0)}(x) \equiv 0$.

On the other hand, we have from (\ref{row:4}) and (\ref{eq:VVV}) that
for $k\geq 2$
\begin{eqnarray*}
G_i^{(k)} &=& G_1^{(1)} +
\sum\limits_{j=i+1}^n \frac{1}{j(j-1)}\left[ W_j^{(k-1)}-
  V_j^{(k-1)}\right] \\
&=&G_i^{(1)}+\frac{1}{i}\left[ V_i^{(k-1)}-\sum _{j=i+1}^n\frac{i}{j(j-1)}
V_j^{(k-1)}\right] ,
\end{eqnarray*}
by (\ref{row:5}).
On letting $i/n \rightarrow x$ as
$n\rightarrow\infty$, we derive from the case $k=m$ that
\[
  G^{(m)}(x) = G^{(1)}(x) + \left(\frac{1}{x}\right) \left[
    V^{(m-1)}(x) - \int\limits_x^1 \frac{x}{y^2} V^{(m-1)}(y)dy
    \right].
\]
Then $G^{(m)}(s_m^*) = 0$ implies that
\[
  G^{(1)}(s_m^*) + \frac{1}{s_m^*} V^{(m-1)}(s_m^*) -
  \int\limits_{s_m^*}^1 \frac{1}{y^2} V^{(m-1)}(y) dy = 0,
\]
or equivalently, from (\ref{row:16}) and $V^{(m-1)}(s_m^*) = q_{m-1}^*$,
\begin{equation}
  \label{row:30}
  \int\limits_{s_m^*}^1 \frac{1}{y^2} V^{(m-1)}(y) dy =  -\frac{(2+\ln s_m^*)\ln s_m^*}{2}
  + \frac{q_{m-1}^*}{s_m^*}.
\end{equation}
Application of (\ref{row:30}) to (\ref{row:29}) provides
$  q_m^* = -s_m^* \ln s_m^* + q_{m-1}^*$, which upon
repetition and use of (\ref{eq:qs}) provides (\ref{row:28}).
\end{Proof}

Numerical values of the first four $q_m^*$ are $q_1^* =
0.2707$, $q_2^* = 0.4725$, $q_3^* = 0.6208$, $q_4^* = 0.7287$. See
Table \ref{tab2} for $q_5^*$ and $q_{10}^*$ ($c=0$).

\begin{remark} The above approach is rather
intuitive. To make the argument more rigorous, we can approximate the difference
equations by
differential equations. This method was suggested by ~\cite{dynyus69:Markov} and
has since been applied
successfully by ~\cite{sza82b}, ~\cite{sucsza02} and ~\cite{yas83}.
~\cite{muc73a,muc73} has developed the idea
for a wider class of optimal stopping problems.

Since
$$V_i^{(k)}=\sum _{j=i+1}^n\frac{i}{j(j-1)}W_j^{(k)}\quad \mbox{ and }\quad
V_{i-1}^{(k)}=\sum _{j=i}^n\frac{i-1}{j(j-1)}W_j^{(k)},$$
we have by subtraction that
$$
V_i^{(k)}-V_{i-1}^{(k)}=\sum _{j=i+1}^n\frac{1}{j(j-1)}W_j^{(k)}
-\frac{1}{i}W_i^{(k)}
=\frac{1}{i}\left[ V_i^{(k)}-W_i^{(k)}\right] .
$$
Also
$U_i^{(k)}=U_i^{(1)}+V_i^{(k-1)}$,
so
$$
W_i^{(k)}=\max \left( U_i^{(k)},V_i^{(k)}\right)
=\max \left( V_i^{(k)},V_i^{(k)}+U_i^{(1)}\right) .
$$
Hence
\begin{eqnarray*}
V_i^{(k)}-V_{i-1}^{(k)}&=&-\frac{1}{i}\left[ U_i^{(1)}+V_i^{(k-1)}-
V_i^{(k)}\right] ^+ \\
&=&-\left[ \sum _{j=i}^n\frac{1}{j}+\frac{1}{i}V_i^{(k-1)}-
\frac{1}{i}V_i^{(k)}\right] ^+ .
\end{eqnarray*}

In the development of the differential equation approach cited above
it can be shown that asymptotically $V_i^{(k)}/n \approx f^{(k)}(i/n)$.
With
$i/n\rightarrow x$ as $n\rightarrow\infty$, we derive
\[
  \frac{d}{dx}f^{(k)}(x) = - \left(\frac{1}{x}\right) \left[-x\ln
    x + f^{(k-1)}(x) - f^{(k)}(x)\right]^+,\ \ 0 \leq x \leq 1
\]
with boundary condition $f^{(k)}(1) = 0$. Here the  nonincreasing
sequence of critical numbers is
\[
  s_k^*:\ f^{(k)}(s_k^*) = f^{(k-1)}(s_k^*) - s_k^*\ln s_k^*,\ \
    k=1,2,\dots \ .
\]
The function $f^{(k)}(\cdot)$ is constant on $[0, s_k^*]$, so the
expected payoff is $q_k^* = f^{(k)}(0) = f^{(k)}(s_k^*)$. For
example, routine algebra yields for $k=1$ and 2 that
\begin{eqnarray*}
  f^{(1)}(x) & = & \begin{cases}
    -s_1^*\ln s_1^*,&       0\leq x\leq s_1^*\\
    \frac{x}{2}(\ln x)^2,&  s_1^*\leq x \leq 1
    \end{cases},\\[0.5cm]
  f^{(2)}(x) & = & \begin{cases}
    -s_1^*\ln s_1^*-s_2^*\ln s_2^*,&  0 \leq x \leq s_2^*\\
    2s_1^*-\frac{2}{3}x+\frac{x}{2}(\ln x)^2,&
                                        s_2^* \leq x \leq s_1^*,\\
    \frac{x}{2}(\ln x)^2-\frac{x}{6}(\ln x)^3,&  s_1^*\leq x \leq
    1,
    \end{cases}
\end{eqnarray*}
where $s_1^* = \exp(-2)$ and $s_2^* = \exp\{-(1+\sqrt{7/3})\}$.
For $ k \geq 3$, we can proceed in similar way.
\end{remark}

\vskip-2cm

\section{\label{sec3}The multiple-choice duration problem with acquisition costs}
In this section, the multiple-choice duration problem is generalized by imposing a constant 
acquisition cost $c=c(n)>0$ each time an object is chosen. The objective of this problem is to 
maximize the expected net payoff, that is, total possession time less the total acquisition cost 
incurred. 
\subsection{The degenerate distribution of the number of objects}
For simplicity we restrict attention to the case $P\{ M=n\} =1$, so that $\pi _i=1$ for 
$1\leq i\leq n$. To avoid triviality we assume $n>1$.

Consider first the one--choice problem. The expected net payoff resulting from a choosing a 
candidate presenting at time $i$ is
$${\mathcal U}_i^{(1)}:=U_i^{(1)}-c,$$
which by (\ref{dtilde}) is given by
\begin{equation}
{\mathcal U}_i^{(1)}=i(H_n-H_{i-1})-c.
\label{eq:6.U}
\end{equation}
We have that
$${\mathcal U}_{i+1}^{(1)}-{\mathcal U}_i^{(1)}=H_n-H_i-1,$$
which is strictly decreasing in $i$ and is negative for $i=n-1$.

Put
\begin{equation}
K(n):=\min \{ i:H_n-H_i\leq 1\} .
\label{KH}
\end{equation}
Then
$${\mathcal U}_1^{(1)}\leq {\mathcal U}_2^{(1)}\leq \ldots \leq {\mathcal U}_{K(n)}^{(1)}\quad
\mbox{ and }\quad {\mathcal U}_{K(n)}^{(1)}\geq {\mathcal U}_{K(n)+1}^{(1)}\geq \ldots
\geq {\mathcal U}_n^{(1)}.$$
If ${\mathcal U}_{K(n)}^{(1)}\leq 0$, it is optimal never to choose a candidate,
so without loss of generality we may assume ${\mathcal U}_{K(n)}^{(1)}>0$,
that is,
\begin{equation}
0< c<U_{K(n)}^{(1)}.
\label{CK}
\end{equation}
 Further,
there exist integers $a=a(n,c)$, $b=b(n,c)$ satisfying
$$1\leq a\leq K(n)\leq b\leq n$$
such that ${\mathcal U}_i^{(1)}\geq 0$ if and only if $a\leq i\leq b$
and to maximize expected total payoff we never choose a candidate
 presenting at time $i$ when $a\leq i\leq b$ fails. Clearly this holds
also in the $m$--choice problem.

We define $\widetilde{\mathcal V}_i^{(1)}$ as the expected net payoff when
we reject a candidate appearing at time $i\leq b$
but select the next candidate (if any) appearing no later than time $b$.
We then have
\begin{eqnarray}
\widetilde{\mathcal V}_i^{(1)}&=&\sum _{j=i+1}^{b}p(i,j){\mathcal U}_j^{(1)}
\nonumber \\
&=&i\left[ \sum _{j=i+1}^{b}\frac{1}{j-1}(H_n-H_{j-1})-c\left( \frac{1}{i}
-\frac{1}{b}\right) \right] \quad \mbox{for $i\leq b$}.
\label{eq:6.V}
\end{eqnarray}

We now turn attention to the $m$--choice problem. For $i\leq b$
we employ the notation ${\mathcal U}_i^{(k)}$, ${\mathcal V}_i^{(k)}$,
$\widetilde{\mathcal V}_i^{(k)}$, ${\mathcal W}_i^{(k)}$ analogously to
$U_i^{(k)}$, $V_i^{(k)}$,
$\widetilde{V}_i^{(k)}$, $W_i^{(k)}$, respectively
and referring to expected net maximal payoff rather than expected total time
of possession of candidates and with choice of second and subsequent
candidates occurring no later than time $b$. If $m>b-a$ it is clearly optimal to simply choose 
every candidate appearing in $I$, so we suppose $m\leq b-a$. The following theorem 
summarizes the optimal strategy for the $m$--choice problem with acquisition cost.

\begin{Theorem}\label{thm41}
For the $m$--choice duration problem with acquisition cost $c$ subject to (\ref{CK}), there exists 
a sequence $(s_1(c), s_2(c),\ldots ,s_m(c))$ of integral critical numbers such that, whenever 
there remain $k$ choices to be made, the optimal strategy selects the first candidate to appear at 
or after time $s_k(c)$ but no later than $b$. Moreover $s_k(c)$ is nonincreasing in $k$ and 
determined by Theorem \ref{thm11} with $N=b$ and 
\begin{equation}
G_i^{(1)}=H_n-H_{i-1}-\sum _{j=i+1}^{b}\frac{1}{j-1}\left( H_n-H_{j-1}\right)- \frac{c}{b}.
\label{exp}
\end{equation}
 Finally, $s_m(c)\geq a$.
\end{Theorem}

\begin{Proof}
There is nonnegative expected payoff from a candidate selected at time $i$
with $a\leq i\leq b$, but not for one selected after time $b$, so it suffices
to establish the result for candidates arriving at times $i\leq b$.

From (\ref{eq:6.U}) and (\ref{eq:6.V}), we can verify that (\ref{exp}) is
equivalent to
$$
G_i^{(1)}=\frac{1}{i}\left[ {\mathcal U}_i^{(1)}-\widetilde{\mathcal V}_i^{(1)}\right].
$$
We derive that for $i<b$,
$$G_{i+1}^{(1)}-G_i^{(1)}=\frac{1}{i}\left[ H_n-H_{i-1}\right] .$$
The right--hand side is strictly decreasing in $i$ for $i<K(n)$ and
nonpositive for $K(n)\leq i<b$. It follows that $(C1)$--$(C3)$ of Theorem \ref{thm11} are 
satisfied provided that $G_{b}^{(1)}\geq 0$.

To see that this requirement is met, observe that ${\mathcal U}_{b}^{(1)}\geq 0$,
that is,
$$
b\left[ H_n-H_{b(n)-1}\right] -c \geq 0.
$$
Hence
$$
G_{b}^{(1)}=H_n-H_{b-1}-\frac{c}{b}\geq 0
$$
as required.

Thus the conditions of Theorem \ref{thm11} are met. Establishing the theorem now
follows closely the rest of the proof of Theorem \ref{thm21}, operating on the interval $[1,b]$ 
instead of $[1,n]$. Since a candidate arriving before time $a$ is never accepted, we have
finally $s_m(c)\geq a$.
\end{Proof}

The case $c=0$ corresponds to the duration problem treated in Corollary \ref{corol11}.
Thus we have $a(n,0)=1$ and $b(n,0)=n$. In Section 8 we shall need to
compare quantities occurring in that context and the present one. Accordingly
we shall where necessary for clarity write the $G_i^{(k)}$ occurring in
this section as $G_i^{(k)}(c)$ and that of Corollary \ref{corol11} as $G_i^{(k)}(0)$,
etc. We shall need the following result.

\begin{Corollary}\label{corol13}
For $c>0$ and $i\leq b(n,c)$, the value of $g_i(c)=G_i^{(1)}(0)-G_i^{(1)}(c)$
is independent of $i$ and so may be written $g(c)$. Further
\begin{equation}
g(c)>c/n.
\label{eq:gd}
\end{equation}
\end{Corollary}
\begin{Proof}
We have for $i\leq b$ that
$$g_i(c)=c/b\ -\sum _{j=b+1}^n\frac{1}{j-1}\left( H_n-H_{j-1}\right) ,$$
which is independent of $i$. Also $j>b$ implies that $j(H_n-H_{j-1})<c$, so
$$\sum _{j=b+1}^n\frac{1}{j-1}(H_n-H_{j-1})<\sum _{j=b+1}^n\frac{c}{j(j-1)}
=c\left( \frac{1}{b}-\frac{1}{n}\right) ,$$
from which (\ref{eq:gd}) follows.
\end{Proof}

\subsection{\label{sec3asym}Asymptotics for the duration problem with acquisition costs}

Observe first that, from (\ref{KH}), $\us{\limn}K(n)/n=e^{-1}$, so the cost
condition (\ref{CK}) is reduced, as $n\rightarrow\infty$, to
\begin{equation}\label{row:42}
  c=\lim_{n\to \infty} c(n)/n\leq e^{-1}.          
\end{equation}
After division by $n$, we may let $i/n\rightarrow x$ as $n\rightarrow\infty$ in
 (\ref{eq:6.U}) to show that $\overline{{\mathcal U}}_i^{(1)}$ approaches
\begin{equation*}
  U^{(1)}(x)=-c+x\int_x^1 \frac{dy}{y}=-c-x \ln x.
\end{equation*}
Let $\beta=\us{\limn}b(n,c)/n$. Then $\beta$ is the unique root $x \in [e^{-1},1)$ of $U^{(1)}(x)=0$ under the cost condition (\ref{row:42}) and satisfies
$$
 -\beta\ln\beta=c.
$$
 For $k\geq 1$, define $s_k^*=s_k^*(c)=\us{\limn}s_k(c)/n$. When we let
 $i/n\rightarrow x$ as $n\rightarrow \infty$ in (\ref{exp}) divided by $n$,
 $G^{(1)}(x)$ approaches the integral
 \begin{equation}
  G^{(1)}(x)=\int_x^1\frac{dy}{y}-\int_x^{\beta}\frac{dy}{y}
\int_y^1\frac{dz}{z}-\frac{c}{\beta}=
  -\frac{(2+\ln x)\ln\,x-(2+\ln\,\beta)\ln \beta}{2}.
 \label{row:44}
 \end{equation}

Thus $s_1^*=\exp\{-(2+\ln \beta)\}$ is obtained as the unique root $x\in(0,\,\beta)$ of
$G^{(1)}(x)=0$. For $0<x<\beta$ and $k\geq 1$, define G$^{(k)}(x)$ recursively by
\begin{equation*}
G^{(k)}(x)=  G^{(1)}(x)+\int_{\max(x,s_{k-1}^*)}^\beta\frac{1}{y}\, 
G^{(k-1)}(y)\,dy,\mbox{ $k\geq2$}  
\end{equation*}
starting with $G^{(1)}(x)$. Then, for G$^{(k-1)}(x)$ and $s_{k-1}^*$ given,
we obtain $s_k^*$ as the unique root $x\in(0,\,s_{k-1}^*)$ of $G^{(k)}(x)=0$, 
or equivalently, from (\ref{row:44}),
\begin{equation*}
s_k^*=\exp\left\{-\left(1+\sqrt{(1+\ln \beta)^2+2\int_{s_{k-1}^*}^\beta
  \frac{G^{(k-1)}(y)}{y}\,dy}\right)\right\}.
\end{equation*}

Similarly to the development of Section 5, we obtain the following result which gives
a generalized version of formula (\ref{row:1}).
\begin{Lemma}\label{lem31}
  Under (\ref{row:42}), $s_k^*$ satisfies the recursion
$$
 s_k^*=\exp\left[-\left\{1+\sqrt{(1+\ln \beta)^2-2\sum\limits_{i=1}^{k-1}
\frac{[(k-i+2)B_{k+1,\,i}+(k-i+1)B_{k+2,\,i}]}{(k-i+2)!}}\right\}\right],
$$
where $B_{k,i}=(\ln\, s_i^*)^{k-i}-(\ln\, \beta)^{k-i}$.
\end{Lemma}

Let $\alpha=1+\ln \beta$. Then from Lemma~\ref{lem31} we can calculate the $s_k^*$ 
successively as
\begin{eqnarray*}
 s_1^*&=&\exp\{-(1+\alpha)\},\\
 s_2^*&=&\exp\left\{-\left(1+\alpha\sqrt{1+\frac{4}{3}\alpha}\right)\right\},\\
 s_3^*&=&\exp\left\{-\left[1+\alpha\sqrt{1+\frac{2}{3}\alpha
          \{1+(1+\frac{4}{3}\alpha)^{3/2}\}}\right]\right\}.
\end{eqnarray*}
For $m\geq1$, let $q_m^*$ be the scaled expected net payoff for the $m$-choice
duration problem when $n$ tends to infinity. Then we have the following result.

\begin{Lemma}\label{lem32}
 Under (\ref{row:42}), we have for $m\geq 1$ that
$$
 q_m^*=-\left( \sum_{k=1}^m s_k^*\ln s_k^*+mc\right).
$$
\end{Lemma}

\begin{Proof} 
Similar to that of  Lemma \ref{lem24}.
\end{Proof}

\begin{table}[bth2]
\renewcommand{\arraystretch}{.8}
 \begin{center}
  \parbox{0.72\textwidth}{
\caption{\label{tab1}The asymptotic critical number $s_m^*$ for some values of
$m$ and $c$.}}
 \end{center}
 $$
\begin{tabular}{l@{\ \ \ \ }c@{\ \ \ \ }c@{\ \ \ \ }c@{\ \ \ \ }c@{\ \ \ \ }c@{\
\ \ \ }
  c@{\ \ \ \ }c@{\ \ \ \ }} \noalign{\vskip1mm}\hline\noalign{\vskip1mm}
$c$& $\beta$& $s_1^*$& $s_2^*$& $s_3^*$& $s_5^*$& $s_{10}^*$&
$s_{\infty}^*(=\beta')$\\
 \noalign{\vskip1mm} \hline\noalign{\vskip1mm}
 0.0& 1.0000&    0.1353& 0.0799& 0.0493& 0.0199& 0.0024& 0.0000\\
 0.1& 0.8942&    0.1513& 0.0990& 0.0698& 0.0416& 0.0281& 0.0280\\
 0.2& 0.7717&    0.1754& 0.1294& 0.1047& 0.0839& 0.0787& 0.0787\\
 0.3& 0.6130&    0.2208& 0.1898& 0.1761& 0.1690& 0.1684&
 0.1684\\[1mm]
 \hline\hline
 \end{tabular}$$
 \end{table}

 \begin{table}[bth2]
\renewcommand{\arraystretch}{.8}
 \begin{center}
  \parbox{0.77\textwidth}{
\caption{\label{tab2}The asymptotic expected net payoff for some values of $m$
and $c$.}}
 \end{center}
 $$
\begin{tabular}{l@{\ \ \ \ }c@{\ \ \ \ }c@{\ \ \ \ }c@{\ \ \ \ }c@{\ \ \ \ }c@{\
\ \ \ }
  c@{\ \ \ \ }} \noalign{\vskip1mm}\hline\noalign{\vskip1mm}
  $c$& $q_1^*$& $q_2^*$& $q_3^*$& $q_5^*$& $q_{10}^*$ &  $q_{\infty}^*$\\
  \noalign{\vskip1mm} \hline\noalign{\vskip1mm}
   0.0&    0.2707& 0.4725& 0.6208& 0.8066& 0.9656& 1.0000\\
   0.1&    0.1858& 0.3147& 0.4005& 0.4871& 0.5195& 0.5197\\
   0.2&    0.1053& 0.1700& 0.2062& 0.2322& 0.2363& 0.2363\\
   0.3&    0.0335& 0.0489& 0.0547& 0.0569& 0.0570& 0.0570\\[1mm]
 \hline\hline
 \end{tabular}$$
 \end{table}

Table \ref{tab1} presents numerical values of $\beta$ and $s_m^*$ for  some
values of $m$ and $c$.
Let $\beta'$ be the unique root  $x\in(0,\,e^{-1}]$ of $-x \ln
x=c$. It is
intuitively clear that, as $m\rightarrow\infty$, $s_m^*$ converges to
$\beta'$ $(=s_\infty^*)$ because
 there is no benefit in choosing a candidate prior to $\beta'$.

Table \ref{tab2} presents numerical values of $q_m^*$ for some values of $m$ and $c$.
It is interesting to compare, for example, $q_1^*=0.0335$ for $c=0.3$ to $q_1^*=0.2707$ for $c=0$,
which implies that we can still gain positive expected payoff even when the acquisition cost is 
larger than the mean maximum payoff attainable when the acquisition cost is zero. This is not a 
contradiction. The stopping region shrinks as $c$ gets large (see Table \ref{tab1}) and positive 
mean payoff is assured  by restricting our choice to a really good object. Table \ref{tab2} 
suggests also that, as $m\rightarrow\infty$, $q_m^*$ converges to a value $q_\infty^*$. This is 
given in the following lemma.

\begin{Lemma}\label{lem33}
 \begin{equation}
  q_\infty^*=(\beta-\beta')\left(1-\frac{c^2}{\beta\beta'}\right)
 \label{row:48}
 \end{equation}
\end{Lemma}
\begin{Proof}{\bf \ref{lem33}} As the arrival times of the $n$
objects, we consider time epochs $1/n,2/n,\ldots ,n/n$ instead of $1,2,\ldots ,n
$. When $n\to \infty$, the transition probability $p(i/n,\,j/n)=i/(j(j-1))$
then converges to the transition density $p(x,\,y)=x/y^2$ as 
$i/n\rightarrow x$, $j/n\rightarrow y$ (see (\ref{row:2})) and the candidates
appear according to a non-homogeneous Poisson process with intensity function
$\lambda(x)=1/x$ from (a), (b), in Section \ref{sec2}. That is, if we let $N(a,\,b)$ denote the 
number of candidates that appear in time interval $(a,\,b)$, then $N(a,\,b)$ becomes a Poisson 
random variable with parameter $\ln (b/a)$ (see Theorem 1 of ~\cite{gilmos66}).

Let $T(x)$ denote the time of the first candidate after time $x$ if there is one
and $1$ if there is not. From the above $T(x)$ has density $ f_{T(x)}(t)=p(x,\,t)=x/t^2$
on the time interval $(x,\,1)$ and probability mass $x$ at $1$. As the number of choices 
$m\rightarrow\infty$, the optimal strategy chooses all the candidates that appear in time interval
$(\beta,\,\beta')$. Thus the total proportional duration $D$ is expressed as

\begin{equation*}
 D=\begin{cases} T(\beta)-T(\beta'), &\text{if }T(\beta')\leq \beta\\
                 0,&\text{if }T(\beta')>\beta
   \end{cases}.
\end{equation*}
It is readily verified that $T(\beta)$ and $T(\beta')$ are independent. Hence by conditioning on
$T(\beta')$, 
\begin{align}
\nonumber \bE[\,D\,]&=\bE[\,T(\beta)-T(\beta')\,|\,T(\beta')\leq \beta\,]
P\{T(\beta')\leq \beta\}\\[2mm]
\nonumber &=\bE[\,T(\beta)\,] P\{T(\beta')\leq
\beta\}-\bE[\,T(\beta')\,|\,T(\beta')\leq \beta\,]
  P\{T(\beta')\leq \beta\} \\
\nonumber &=\left\{\int_{\beta}^1 t\,f_{T(\beta)}(t)dt+\beta\right\}
   \left\{\int_{\beta'}^{\beta} f_{T(\beta')}(t)dt\right\}
   - \int_{\beta'}^\beta t\,f_{T(\beta')}(t)dt\\
\nonumber
&=(c+\beta)\left(1-\frac{\beta'}{\beta}\right)
-c\left(1-\frac{\beta'}{\beta}\right)\\
\nonumber &=\beta-\beta'.
\end{align}
Thus the expected net payoff $q_\infty^*$ is
\begin{equation*}
\bE[\,D-c\,N(\beta,\beta')\,]=(\beta-\beta')-c \ln \left(\frac{\beta}
{\beta'}\right) ,
\end{equation*}
which yields (\ref{row:48}).
\end{Proof}


\section{\label{sec4}Duration problem with replacement costs }
In this section a constant cost $d=d(n)>0$ is incurred each time there is replacement, whether or 
not the new candidate is the one to end the candidature of the previously chosen candidate.
For simplicity we consider only the case where $M\sim M_s(n)$ and ignore acquisition costs. The 
objective is to maximize the expected net payoff, that is, the total time of possession of a 
relatively best object less any replacement costs incurred. The multiple-choice duration problem 
with a replacement cost may be considered as a marriage and divorce problem, interpreting the 
replacement cost as alimony.

\subsection{The degenerate distribution of the number of objects}
We treat the $m$-choice duration problem with replacement cost $d>0$. In the $m$-choice problem we 
are allowed to replace objects up to $m-1$ times, $m\geq 2$. We define the state of the process
as in Section 3 and ${\mathcal W}_i^{(k)}$, ${\mathcal U}_i^{(k)}$
and ${\mathcal V}_i^{(k)}$ similarly to in Section 6.

Consider a candidate other than the first arriving at time $i$. As in
Section 6, we may argue that such a candidate is never chosen unless
\begin{equation}
d<U_{K(n)}.
\label{dK}
\end{equation}
Further, $U_i^{(1)}-d\geq 0$ if and only if $a(n,d)\leq i\leq b(n,d)$.

Once the  first choice is made, the problem reduces to the $(m-1)$-choice problem with an 
acquisition cost $d$. Thus the main concern is to determine when to make the first choice. The 
optimal strategy can be summarized as follows.

\begin{Theorem}\label{thm51}
     For the $m$-choice duration problem with replacement cost condition (\ref{dK}), there
     exists a sequence $(s_1(d),s_2(d),\ldots,s_{m-1}(d),t_m(d))$ of integral critical numbers   
     such that the optimal strategy first selects the first candidate (if any) to appear at or 
     after time $t_m(d)$. Thereafter it replaces each previously chosen object with the first new 
     candidate (if any) that appears at or after time $s_k(d)$ but no later than $b(n,d)$ if $k$ 
     more replacements are available ($1\leq k\leq m-1$), where $b(n,d)=\max\{i:U_i^{(1)}\geq d\}$.

     Each $s_k(d)$ is as in Theorem \ref{thm41} while $t_m(d)\leq s_{m}(d)$ and is determined by
\begin{equation}
      t_m(d)=\min\left\{i\leq b(n,d):G_i^{(m)}(d)+g(d)\geq 0 \right\}.
\label{row:50B}
\end{equation}
\end{Theorem}

\begin{Proof}
The part of the result relating to choices when fewer than $m$ replacements are to be made is 
immediate from Theorem \ref{thm51}, so it remains to address the first choice of a candidate.

As before
$$
{\mathcal W}_j^{(m-1)}=\begin{cases} {\mathcal V}_j^{(m-1)}, &\text{if }
j< s_{m-1}(d)\\
{\mathcal U}_j^{(m-1)}, &\text{if }
j\geq s_{m-1}(d)\\
   \end{cases}
$$
and
$$
{\mathcal V}_j^{(m-1)}=\widetilde{{\mathcal V}}_j^{(m-1)}\quad
\mbox{ if }j\geq s_{m-1}(d)-1.
$$
The principle of optimality provides
$$
{\mathcal U}_i^{(m)}=U_i^{(1)}+\sum _{j=i+1}^{b(n,d)}p(i,j){\mathcal W}_j^{(m-1)},
$$
$$
\widetilde{{\mathcal V}}_i^{(m)}=\sum _{j=i+1}^{b(n,d)}p(i,j)\left[ U_j^{(1)}+
{\mathcal V}_j^{(m-1)} \right]
=\widetilde{V}_i^{(1)}+\sum _{j=i+1}^{b(n,d)}p(i,j){\mathcal V}_j^{(m-1)},
$$
so that
\begin{eqnarray*}
{\mathcal U}_i^{(m)}-\widetilde{{\mathcal V}}_i^{(m)}&=&
U_i^{(1)}-\widetilde{U}_i^{(1)}+\sum _{j=\max (i+1,s_{m-1}(d))}^{b(n,d)}p(i,j)
\left[ {\mathcal U}_i^{(m)}-\widetilde{{\mathcal V}}_i^{(m)}\right] \\
&=&iG_i^{(1)}(0)-iG_i^{(1)}(d)+iG_i^{(m)}(d).
\end{eqnarray*}
Thus
$$\frac{{\mathcal U}_i^{(m)}-\widetilde{{\mathcal V}}_i^{(m)}}{i}
=G_i^{(m)}(d)+g(d).$$
Since by Corollary \ref{corol13} $g(d)>0$, (\ref{row:50B}) implies that $t_m(d)\leq
s_m(d)$. Further, ${\mathcal U}_i^{(m)}\geq \widetilde{{\mathcal V}}_i^{(m)} $
if and only if $i\geq t_m(d)$. Thus if the choice of a candidate appearing at
time $i\leq b(n,d)$ is preferable to that of the next candidate (if any) before
time $b(n,d)$, then it is preferable to the choice of any subsequent candidate.
This concludes the proof.
\end{Proof}

\subsection{\label{sec4asym}Asymptotics for the duration problem with replacement costs}

As $n\rightarrow\infty$, the cost condition (\ref{dK}) is reduced to
 \begin{equation}
  d=\lim _{n\to \infty}d(n)/n\leq e^{-1}.
 \label{row:51}
 \end{equation}
Let $\delta=\us{\limn}b/n$. Then under condition (\ref{row:51}),
$\delta$ is the unique root
$x\in[e^{-1},\,1)$ of $-x\ln x=d$. We have the following result
concerning the limiting values
$s_k^*=\us{\limn}s_k/n$ ($k\geq 1$) and $t_m^*=\us{\limn}t_m/n$.

 \begin{Lemma}
  Under (\ref{row:51}), $t_m^*$ may be expressed in terms of $s_m^*$ as
 \begin{equation}
t_m^*=\exp\left[-\left\{1+\sqrt{(1+\ln s_m^*)^2-(2+\ln \delta)\ln \delta}
        \right\}\right],
 \label{row:52}
 \end{equation}
 where $s_k^*$ ($1\leq k\leq m$) satisfies the recursion
 \begin{equation}
  s_k^*= \exp\left[-\left\{1+\sqrt{(1+\ln \delta)^2-2
\sum\limits_{i=1}^{k-1}\frac{[(k-i+2)B_{k+1,\,i}+(k-i+1)B_{k+2,\,i}]}{(k-i+2)!}}
        \right\}\right],
 \label{row:53}
 \end{equation}
 with $B_{k,\,i}=(\ln s_i^*)^{k-i}-(\ln \delta)^{k-i}$.
 \end{Lemma}

 \begin{Proof} Equation (\ref{row:53}) is evident from Lemma~\ref{lem31},
while (\ref{row:52}) is immediate from (\ref{row:50B}).
 \end{Proof}

 Let $\lambda=1+\ln \delta$. Then from (\ref{row:52}) and
 (\ref{row:53}) we have
 \begin{eqnarray*}
   t_2^*&=&\exp\left\{-\left(1+\sqrt{1+\frac{4}{3}\lambda^3}\right)\right\},\\
   t_3^*&=&\exp\left\{-\left[1+\sqrt{1+\frac{2}{3}\lambda^3
    \left\{1+\left(1+\frac{4}{3}\lambda\right)^{3/2}\right\}} \right]\right\}. \\
 \end{eqnarray*}

 For $m\geq 2$, let $r_m^*$ be the expected net payoff for the
$m$-choice duration problem when $n$ tends to infinity. Then we have
the following.

 \begin{Lemma}\label{lem42}

\begin{description}
\item[$(i)$] If $d>e^{-1}$, then  $r_m^*=2e^{-2}$.
\item[$(ii)$] If $d\leq e^{-1}$, then $r_m^*=-\left[ \sum_{k=1}^{m-1} s_k^*\ln
s_k^*+t_m^*\ln t_m^*+(m-1)d\right]$.
\end{description}
\end{Lemma}

 \begin{Proof} The proof is omitted.
 \end{Proof}

 \begin{table}[bth2]
  \renewcommand{\arraystretch}{.8}
 \begin{center}
  \parbox{0.72\textwidth}{
\caption{\label{tab3}The asymptotic critical number $t_m^*$ for some values of
$m$ and $d$.}}
 \end{center}
 $$
\begin{tabular}{l@{\ \ \ \ }c@{\ \ \ \ }c@{\ \ \ \ }c@{\ \ \ \ }c@{\ \ \ \ }c@{\
\ \ \ }
  c@{\ \ \ \ }c@{\ \ \ \ }}\noalign{\vskip1mm} \hline\noalign{\vskip1mm}
  $d$& $t_2^*$& $t_3^*$& $t_5^*$& $t_{10}^*$& $t_{\infty}^*$\\
  \noalign{\vskip1mm} \hline\noalign{\vskip1mm}
  0.1&    0.0916& 0.0656& 0.0397& 0.0270& 0.0268&\\
  0.2&    0.1063& 0.0885& 0.0725& 0.0684& 0.0684&\\
  0.3&    0.1243& 0.1186& 0.1154& 0.1151& 0.1151&\\ [1mm]
 \hline\hline
 \end{tabular}$$
 \end{table}

 \begin{table}[bth2]
\renewcommand{\arraystretch}{.8}
 \begin{center}
  \parbox{0.72\textwidth}{
\caption{\label{tab4}The asymptotic expected net payoff $r_m^*$ for some values
of $m$ and $d$.}}
    \end{center}
$$\begin{tabular}{l@{\ \ \ \ }c@{\ \ \ \ }c@{\ \ \ \ }c@{\ \ \ \ }c@{\ \ \ \
}c@{\ \ \ \ }
  c@{\ \ \ \ }c@{\ \ \ \ }}\noalign{\vskip1mm} \hline\noalign{\vskip1mm}
  $d$& $r_2^*$& $r_3^*$& $r_5^*$& $r_{10}^*$& $r_{\infty}^*$\\
  \noalign{\vskip1mm} \hline\noalign{\vskip1mm}
  $d$& $r_2^*$& $r_3^*$& $r_5^*$& $r_{10}^*$& $r_{\infty}^*$\\ \noalign{\vskip2mm}
  \hline\noalign{\vskip2mm}
  0.1&    0.4047& 0.4934& 0.5828& 0.6166& 0.6168\\
  0.2&    0.3435& 0.3845& 0.4146& 0.4198& 0.4198\\
  0.3&    0.2927& 0.3017& 0.3056& 0.3059& 0.3059\\[1mm]
\hline\hline
\end{tabular}$$
\end{table}

 Tables \ref{tab3} and \ref{tab4} give numerical values of $t_m^*$ and $r_m^*$
for some values of $m$ and $d$ respectively. The values of $s_m^*$ are given in
Table \ref{tab1}
if $c$ is interpreted as $d$. Tables \ref{tab3} and \ref{tab4} suggest that, as
$m\rightarrow\infty$,
$t_m^*$ and $r_m^*$ converge to limits $t_\infty^*$ and $r_\infty^*$
respectively. The following
 lemma specifies these.

 \begin{Lemma}
Let $\delta'$ be the unique root $x\in(0,\,e^{-1}]$ of $-x\ln x=
d$ for $d\leq e^{-1}$. Then
 \begin{equation}
r_\infty^*=(\delta-\delta')\left(1-\frac{d^2}{\delta\delta'}\right)
-t_\infty^*\ln t_\infty^*,
 \label{row:54}
 \end{equation}
 where
 \begin{equation}
t_\infty^*= \exp
\left\{-\left[1+\sqrt{1-2d\left(\frac{\delta-\delta'}{\delta\delta'}\right)
+(\delta-\delta')(\delta+\delta')\left(\frac{d}{\delta\delta'}\right)^2}
\right]\right\} .
 \label{row:55}
 \end{equation}
 \end{Lemma}

 \begin{Proof} Relation (\ref{row:55}) is immediate from (\ref{row:52}), while
(\ref{row:54}) is immediate from Lemmas
 \ref{lem32}, \ref{lem33} and \ref{lem42} (ii).
 \end{Proof}

\vskip-3cm

\section{\label{sec5}Final remarks}

The closely-related multiple-choice secretary problems have been considered by
~\cite{nik76:minimal,nik77:odnom},
~\cite{tam79:double}, ~\cite{mor84:rand}, ~\cite{sta85}, ~\cite{sak87},
~\cite{wil91:assignment},
~\cite{ano89:three} and others. ~\cite{pre89:choice} gives some interesting
results and a review of the literature. There are also
results for Dynkin game models of
the secretary problem when one player has the opportunity to stop and accept a
 candidate more than once (see
~\cite{ksz2000:more} and ~\cite{yassza02:dynkin}).

Multiple-stopping models have recently been applied as a modelling tool for technical and economic
phenomena. ~\cite{sza96:dis} has investigated the double-disorder problem for discrete-time Markov 
processes. ~\cite{assgolsam04:two} and ~\cite{kuhrus02:two} have considered the asymptotic 
properties of the double stopping procedure for iid random variables with known distribution.
~\cite{assgolsam02:several,asssam00:mortal} studied prophet inequalities in the case when the 
mortal has several choices.

The model presented for the multiple-exchange duration problem is important for both  applications and theoretical investigation of the optimization techniques in a stochastic environment. There is opportunity for further extensions of the model in many directions (cf. ~\cite{maztam03:dur,maztam06:dur}).

Some analytical aids exist for numerical and theoretical work, though we have not needed to invoke 
them in this study. We note in particular that the harmonic number $H_n$ can be expressed 
analytically as $H_n=\gamma+\psi _0(n+1)$, where $\gamma\cong 0.577216$ is the Euler-Mascheroni
constant (see \cite{graryz00:functions}) and the digamma function $\psi _0$ is defined by
$$
\psi_0(z)=\frac{d}{d\!z}\ln\Gamma(z)= \frac{\Gamma^{'}(z)}{\Gamma(z)}
$$
(see \cite{abrste72:functions} p. 260).
\vskip-1cm

\vskip-1cm

\section*{\bf Acknowledgement}

\vskip-0.75cm

We should like to thank Professor Masami Yasuda for helpful suggestions. We also thank Masuyo Kawai
for help with the numerical evaluations. The first author was supported by Grant-in-Aid for 
Scientific Research (c) 10680439. Finally, thank are due to an anonymous referee for comments on 
an earlier draft.

\vskip-2cm


\begin{thebibliography}{49}
\expandafter\ifx\csname natexlab\endcsname\relax\def\natexlab#1{#1}\fi
\expandafter\ifx\csname url\endcsname\relax
  \def\url#1{\texttt{#1}}\fi
\expandafter\ifx\csname urlprefix\endcsname\relax\def\urlprefix{URL }\fi

\bibitem[{Abramowitz and Stegun(2000)}]{abrste72:functions}
Abramowitz, M., Stegun, I.~A. (Eds.), 2000. Handbook of Mathematical Functions
  with Formulas, Graphs, and Mathematical Tables. Dover, New York.

\bibitem[{Ano(1989)}]{ano89:three}
Ano, K., 1989. {O}ptimal selection problem with three stops. J. Oper. Res. Soc.
  Japan 32, 491--504.

\bibitem[{Assaf et~al.(2002)Assaf, Goldstein, and
  Samuel-Cahn}]{assgolsam02:several}
Assaf, D., Goldstein, L., Samuel-Cahn, E., 2002. {Ratio prophet inequalities
  when the mortal has several choices.} Ann. Appl. Probab. 12~(3), 972--984.

\bibitem[{Assaf et~al.(2004)Assaf, Goldstein, and
  Samuel-Cahn}]{assgolsam04:two}
Assaf, D., Goldstein, L., Samuel-Cahn, E., 2004. {Two-choice optimal stopping.}
  Adv. Appl. Probab. 36~(4), 1116--1147.

\bibitem[{Assaf and Samuel-Cahn(2000)}]{asssam00:mortal}
Assaf, D., Samuel-Cahn, E., 2000. {Simple ratio prophet inequalities for a
  mortal with multiple choices.} J. Appl. Probab. 37~(4), 1084--1091.

\bibitem[{Bearden(2006)}]{bea06:cardinal}
Bearden, J.~N., 2006. A new secretary problem with rank-based selection and
  cardinal payoffs. J. Math. Psychology 50, 58 -- 59.

\bibitem[{Dynkin and Yushkevich(1969)}]{dynyus69:Markov}
Dynkin, E.~B., Yushkevich, A.~A., 1969. {M}arkov {P}rocess, {T}heorems and
  {P}roblems. Plenum Press, New York.

\bibitem[{Ehjdukyavichyus(1979)}]{eid79}
Ehjdukyavichyus, R., 1979. {O}ptimalna ostanovka markovskoj cepi dvumia
  momentami ostanovki. Litov. Mat. Sb. 13, 181--183.

\bibitem[{Ferguson(1989)}]{fer89:who}
Ferguson, T.~S., 1989. {W}ho solved the secretary problem? Statist. Sci. 4,
  282--296.

\bibitem[{Ferguson et~al.(1992)Ferguson, Hardwick, and
  Tamaki}]{ferhartam92:own}
Ferguson, T.~S., Hardwick, J.~P., Tamaki, M., 1992. {M}aximizing the duration
  of owning a relatively best object. Contemp. Math. 125, 37--57.

\bibitem[{Gilbert and Mosteller(1966)}]{gilmos66}
Gilbert, J.~P., Mosteller, F., 1966. {R}ecognizing the maximum of a sequence.
  J. Am. Stat. Assoc. 61, 35--73.

\bibitem[{Gnedin(2004)}]{gne04:planar}
Gnedin, A.~V., 2004. Best choice from the planar poisson process. Stoch. Proc.
  Appl. 111, 317--354.

\bibitem[{Gnedin(2005)}]{gne05:objectives}
Gnedin, A.~V., 2005. {O}bjectives in the best-choice problems. Sequential
  Analysis 24, 1--11.

\bibitem[{Gradshteyn and Ryzhik(2000)}]{graryz00:functions}
Gradshteyn, I.~S., Ryzhik, I.~M., 2000. Tables of Integrals, Series, and
  Products, 6th Edition. Academic Pres, San Diego, CA.

\bibitem[{Haggstrom(1967)}]{hag67}
Haggstrom, G.~W., 1967. {O}ptimal sequential procedures when more then one stop
  is required. Ann. Math. Stat. 38, 1618--1626.

\bibitem[{K\"{u}hne and R\"{u}schendorf(2002)}]{kuhrus02:two}
K\"{u}hne, R., R\"{u}schendorf, L., 2002. On optimal two-stopping problems. In:
  Berkes, I., Cs\'aki, E., R\'ev\'esz, P. (Eds.), {Limit theorems in
  probability and statistics}. Vol.~II. J\'{a}nos Bolyai Mathematical Society,
  Budapest, pp. 261--271.

\bibitem[{Lehtinen(1993)}]{leh93:unknown}
Lehtinen, A., 1993. {T}he best-choice problem with an unknown number of
  objects. Z. Oper.Res. 37, 97--106.

\bibitem[{Mazalov and Tamaki(2003)}]{maztam03:dur}
Mazalov, V.~V., Tamaki, M., 2003. Explicit solutions to the duration problem.
  Aichi Keiei Ronsyu 147, 69--92.

\bibitem[{Mazalov and Tamaki(2006)}]{maztam06:dur}
Mazalov, V.~V., Tamaki, M., 2006. An explicit formula for the optimal gain in
  the full-information problem of owning a relatively best object. J. Appl.
  Probab. 43~(1), 87--101.

\bibitem[{M\'ori(1984)}]{mor84:rand}
M\'ori, T.~F., 1984. {T}he random secretary problem with multiple choice.
  Annales Univ. Sci. Budapestinensis de Rolando Eotvos Nominatae V, 91--102.

\bibitem[{Mucci(1973{\natexlab{a}})}]{muc73a}
Mucci, A.~G., 1973{\natexlab{a}}. {D}ifferential equations and optimal choice
  problem. Ann.Stat. 1, 104--113.

\bibitem[{Mucci(1973{\natexlab{b}})}]{muc73}
Mucci, A.~G., 1973{\natexlab{b}}. {O}n a class of secretary problems. Ann.
  Probab. 1, 417--427.

\bibitem[{Nikolaev(1976)}]{nik76:minimal}
Nikolaev, M.~L., 1976. On the selection of two objects with minimal sum rank.
  Izv. Vyssh. Uchebn. Zaved., Mat. 3(166), 33--42, {Z}adacha vybora dvokh
  ob"ektov z minimal'nom sumarnym rangom.

\bibitem[{Nikolaev(1977)}]{nik77:odnom}
Nikolaev, M.~L., 1977. {On a generalization of the best choice problem.} Theory
  Probab. Appl. 22, 187--190.

\bibitem[{Nikolaev(1979)}]{nik79}
Nikolaev, M.~L., 1979. {O}bobshchennyje posledovatelnyje procedury. Litov. Mat.
  Sb. 191, 35--44.

\bibitem[{Nikolaev(1998)}]{nik98:multstopping}
Nikolaev, M.~L., 1998. {Optimal multi-stopping rules.} Obozr. Prikl. Prom. Mat.
  5~(2), 309--348.

\bibitem[{Petruccelli(1983)}]{pet83}
Petruccelli, J.~D., 1983. {O}n the best-choice problem when the number of
  observation is random. J. Appl. Probab. 20, 165--171.

\bibitem[{Porosi\'nski(1987)}]{por87:full}
Porosi\'nski, Z., 1987. {T}he full-information best choice problem with a
  random number of observations. Stochastic Processes and their Applications
  North-Holland 24, 293--307.

\bibitem[{Porosinski(2002)}]{por02:similar}
Porosinski, Z., 2002. {On best choice problems having similar solutions.} Stat.
  Probab. Lett. 56~(3), 321--327.

\bibitem[{Preater(1994)}]{pre89:choice}
Preater, J., 1994. {O}n multiple choice secretary problem. Math. Oper. Res.
  19~(3), 597--602.

\bibitem[{Presman and Sonin(1972)}]{preson72:eng}
Presman, E.~L., Sonin, I.~M., 1972. The best choice problem for a random number
  of objects. Theory Prob.Appl. 17, 657--668.

\bibitem[{Sakaguchi(1978)}]{sak78d}
Sakaguchi, M., 1978. {D}owry problem and {O}{L}{A} policies. Rep. Stat. Appl.
  Res. Union Jap. Sci. Eng. JUSE 25, 124--128.

\bibitem[{Sakaguchi(1987)}]{sak87}
Sakaguchi, M., 1987. {G}eneralized secretary problems with three stops. Math.
  Japonica 32, 105--122.

\bibitem[{Samuels(1991)}]{sam91:secretary}
Samuels, S.~M., 1991. {S}ecretary problems. In: Ghosh, B.~K., Sen, P.~K.
  (Eds.), Handbook of Sequential Analysis. Marcel Decker, New York, pp.
  381--405.

\bibitem[{Samuels(2004)}]{sam04:why}
Samuels, S.~M., 2004. {Why do these quite different best-choice problems have
  the same solutions?.} Adv. Appl. Probab. 36~(2), 398--416.

\bibitem[{Seale and Rapoport(1997)}]{searap97:exper}
Seale, D.~A., Rapoport, A., 1997. Sequential decision making with relative
  ranks: {A}n experimental investigation of the ''secretary problem''.
  Organizational Behaviour and Human Decision Processes 69, 221--236.

\bibitem[{Seale and Rapoport(2000)}]{searap00:unknown}
Seale, D.~A., Rapoport, A., 2000. Optimal stopping behavior with relative
  ranks: {T}he secretary problem with unknown population size. J. Behavioral
  Decision Making 13, 391--411.

\bibitem[{Stadje(1985)}]{sta85}
Stadje, W., 1985. {O}n multiple stopping rules. Optimization 16, 401--418.

\bibitem[{Suchwa\l{}ko and Szajowski(2002)}]{sucsza02}
Suchwa\l{}ko, A., Szajowski, K., 2002. Non standard, no information secretary
  problems. Sci. Math. Japonicae 56, 443 -- 456.

\bibitem[{Szajowski(1982)}]{sza82b}
Szajowski, K., 1982. Optimal choice problem of $a$-th object. Matem.Stos. 19,
  51--65, in Polish.

\bibitem[{Szajowski(1996)}]{sza96:dis}
Szajowski, K., 1996. A two-disorder detection problem. Applicationes
  Mathematicae 24~(2), 231--241.

\bibitem[{Szajowski(2002)}]{ksz2000:more}
Szajowski, K., 2002. {O}n stopping games when more than one stop is possible.
  In: Kolchin, V.~F., Kozlov, V.~Y., Mazalov, V.~V., Pavlov, Y.~L., Prokhorov,
  Y.~V. (Eds.), Probability Methods in Discrete Mathematics, Proceedings of the
  Fifth International Petrozavodsk Conference, May 2000. International Science
  Publishers, pp. 57--72.

\bibitem[{Szajowski(2006)}]{sza06:risk}
Szajowski, K., 2006. A rank-based selection with cardinal payoffs and a cost of
  choice. Preprint I-18/2006, Instytute of Mathematics and Computer Science,
  Wybrze\.ze Wyspia\,nskiego 27, 50-370 Wroc{\l}aw,
  \href{http://neyman.im.pwr.wroc.pl/\~szajow/publ2002/pdf/RankStop06.pdf}
  {http://neyman.im.pwr.wroc.pl/\~{}szajow/publ2002/pdf/RankStop06.pdf}.

\bibitem[{Tamaki(1979{\natexlab{a}})}]{tam79:double}
Tamaki, M., 1979{\natexlab{a}}. {A} secretary problem with double choice. J.
  Oper. Res. Soc. Jap. 22~(4), 257--265.

\bibitem[{Tamaki(1979{\natexlab{b}})}]{tam79:ola}
Tamaki, M., 1979{\natexlab{b}}. {O}{L}{A} policy and the best choice problem
  with random number of objects. Math. Japonica 24~(4), 451--457.

\bibitem[{Tamaki et~al.(1998)Tamaki, Pearce, and
  Szajowski}]{tampeasza98:duration}
Tamaki, M., Pearce, C.~E., Szajowski, K., 1998. {Multiple choice problems
  related to the duration of the secretary problem.} RIMS Kokyuroku 1068,
  75--86.

\bibitem[{Wilson(1991)}]{wil91:assignment}
Wilson, J.~G., 1991. {O}ptimal choice and assignment of the best $m$ of $n$
  randomly arriving items. Stoch. Proc. Appl. 39, 325--343.

\bibitem[{Yasuda(1983)}]{yas83}
Yasuda, M., 1983. {O}n a stopping problem involving refusal and forced
  stopping. J. Appl. Probab. 20, 71--81.

\bibitem[{Yasuda and Szajowski(2002)}]{yassza02:dynkin}
Yasuda, M., Szajowski, K., 2002. {D}ynkin games and its extension to a multiple
  stopping model. Bulletin of the Japan Society for Industrial Mathematics
  12~(3), 17--28, in Japanese.

\end{thebibliography}

\normalsize

\def\polhk#1{\setbox0=\hbox{#1}{\ooalign{\hidewidth
  \lower1.5ex\hbox{`}\hidewidth\crcr\unhbox0}}}

\end{document}